\documentclass[letterpaper, 10 pt, conference]{ieeeconf}  

\IEEEoverridecommandlockouts                              

\usepackage{graphics,color} 
\usepackage{hyperref}
\usepackage{epsfig} 
\usepackage{amsmath,mathtools,booktabs}
\usepackage[normalem]{ulem}

\usepackage{cases}
\usepackage{cite}
\usepackage[mathscr]{euscript}
\newtheorem{assumption}{Assumption}
 \let\mathscr\relax
\usepackage[scr]{rsfso}

\usepackage{amssymb}
\allowdisplaybreaks
\usepackage{multirow}

\newtheorem{proposition}{Proposition}

\title{\LARGE \bf 
Robust Continuous-Time Generation Scheduling under \\
Power Demand Uncertainty:
An Affine Decision Rule Approach}

\author{Youngchae Cho \and  Insoon Yang \and Takayuki Ishizaki
\thanks{
This paper is based on results obtained from a project, JPNP24007, commissioned by the New Energy and Industrial Technology Development Organization (NEDO). 
This work was supported by the Information and Communications Technology Planning and Evaluation (IITP) grants funded by MSIT No. 2022-0-00124, No. 2022-0-00480 and No. RS-2021-II211343, Artificial Intelligence Graduate School Program (Seoul National University).}%
\thanks{
Y. Cho and T. Ishizaki are with the Department of Systems and Control Engineering, Institute of Science Tokyo, Tokyo, Japan. 
{\tt\small \{cho,ishizaki\}@lim.sc.e.titech.ac.jp}.
I. Yang is with the Department of Electrical and Computer Engineering, ASRI, Seoul National University, Seoul, South Korea.
{\tt\small insoonyang@snu.ac.kr}.}
}

\begin{document}

\maketitle
\thispagestyle{empty}
\pagestyle{empty}{
}
\begin{abstract}
Most existing generation scheduling models for power systems under demand uncertainty rely on energy-based formulations with a finite number of time periods, which may fail to ensure that power supply and demand are balanced continuously over time. To address this issue, we propose a robust generation scheduling model in a continuous-time framework, employing a decision rule approach. First, for a given set of demand trajectories, we formulate a general robust generation scheduling problem to determine a decision rule that maps these demand trajectories and time points to the power outputs of generators. Subsequently, we derive a surrogate of it as our model by carefully designing a class of decision rules that are affine in the current demand, with coefficients invariant over time and constant terms that are continuous piecewise affine functions of time.
As a result, our model can be recast as a finite-dimensional linear program to determine the coefficients and the function values of the constant terms at each breakpoint, solvable via the cutting-plane method. Our model is non-anticipative unlike most existing continuous-time models, which use Bernstein polynomials, making it more practical. We also provide illustrative numerical examples. 
\end{abstract}

\section{Introduction}
Efficiently balancing power supply and demand  is a fundamental objective of generation scheduling problems, such as unit commitment \cite{padhy2004unit} and economic dispatch \cite{kunya2023review}. 
With power demand becoming increasingly unpredictable over the past few decades due to the integration of renewables, stochastic optimization has been used to develop generation scheduling models that address demand uncertainty endogenously \cite{qiu2022application}.\footnote{In this paper, ``demand" refers to net power demand, defined as power demand minus renewable generation, unless stated otherwise.} 
 
Most of the existing generation scheduling models that account for demand uncertainty, e.g., \cite{tang2020reserve} and \cite{yildiran2023robust}, are formulated in discrete-time settings. These models typically consider a finite number of time periods, treating the energy demand for each period as uncertain. As a result, they can only ensure that power supply and demand are balanced on average within each time period. However, power supply and demand must be balanced continuously over time. Moreover, as noted by \cite{guan2000energy}, a discrete-time generation scheduling model may not guarantee the existence of a feasible power output trajectory that complies with the ramp-rate limits of generators, even when energy demand is known in advance. Simply augmenting the temporal granularity does not necessarily resolve these issues and may increase computational burden \cite{pandvzic2014effect}. 
Thus, it is essential to formulate a continuous-time generation scheduling model that is capable of addressing demand uncertainty and producing a feasible power output trajectory aligned with the actual demand trajectory.

The existing continuous-time generation scheduling models under demand uncertainty include 
\cite{
hreinsson2019continuous, le2022aggregation, khatami2019flexibility, khatami2019stochastic,nikoobakht2020continuous}, where multiple demand trajectories are explicitly considered. \cite{hreinsson2019continuous} and \cite{le2022aggregation} address power systems with energy storage systems and heating, ventilation, and air conditioning (HVAC) systems, respectively. \cite{khatami2019flexibility} determines the reserve capacity and deployment trajectories of generators with different time fidelities. This approach is extended by \cite{khatami2019stochastic} to incorporate energy storage devices. \cite{nikoobakht2020continuous} considers integrated electricity and natural gas transmission systems with demand response programs. To capture demand uncertainty, these models utilize finitely many representative scenarios \cite{le2022aggregation,khatami2019flexibility,khatami2019stochastic}, a scenario tree \cite{hreinsson2019continuous}, and the information-gap decision theory \cite{nikoobakht2020continuous}.

Notably, \cite{le2022aggregation,khatami2019flexibility,khatami2019stochastic,hreinsson2019continuous,nikoobakht2020continuous}
commonly leverage Bernstein polynomials to model both demand trajectories and the corresponding power output trajectories of generators. 
In the Bernstein polynomial approach, a continuous function on the unit interval $\left[0,1\right]$ is approximated by a linear combination of a predetermined number of Bernstein basis polynomials, the coefficients of which are referred to as Bernstein coefficients. 
This method is particularly useful because it provides a constructive proof of the Weierstrass approximation theorem, 
which states that any continuous function on a closed interval can be uniformly approximated by polynomials \cite{levasseur1984probabilistic}.
Furthermore, as highlighted in \cite{parvania2015unit}, Bernstein polynomials offer several additional modeling advantages, such as a straightforward representation of derivatives and the well-known ``convex hull property" \cite{dierckx1995curve}. Studies on continuous-time generation scheduling under uncertainty that utilize Bernstein polynomials also include 
\cite{nikoobakht2019integrated, nourollahi2023continuous, xu2023probabilistic, liu2024continuous}.

However, the Bernstein polynomial approach inherently violates the non-anticipativity constraint, which is a modeling principle asserting that decisions must be made solely based on the information available up to the current time point \cite{shapiro2009lectures}.  This occurs because, for any time interval, the Bernstein coefficients of a demand trajectory is required to obtain those of a power output trajectory. 
Since a demand trajectory (and its Bernstein coefficients) over a non-degenerate closed time interval is not fully known at any single time point except for the right endpoint, non-anticipativity is not satisfied.

To resolve this issue, we propose a robust continuous-time generation scheduling model under demand uncertainty that is non-anticipative, 
unlike those in  \cite{le2022aggregation,khatami2019flexibility,khatami2019stochastic,hreinsson2019continuous,nikoobakht2020continuous}. 
To incorporate non-anticipativity, we adopt a decision rule approach \cite{georghiou2019decision}.

To that end, we first present a general generation scheduling problem aimed at determining a decision rule, defined as a function that maps the demand trajectories of loads and time points to the power outputs of generators. Its objective is to minimize the worst-case total generation cost while maintaining the power supply-demand balance at any time point throughout the planning horizon for any demand trajectory from a given set. The decision rule is explicitly constrained by both continuity and non-anticipativity.

Recognizing that this optimization problem is insufficiently specific to apply a standard algorithm, we derive a surrogate model through two key modifications.

First, we replace the set of demand trajectories for each load with a superset defined by two continuous piecewise affine functions of time. These functions correspond to the 
pointwise maximum and minimum demand trajectories, respectively.

Second, we restrict the decision rule to a carefully designed class that incorporates affine decision rules. Affine decision rules are a widely used solution method for operational planning and control under uncertainty in discrete-time settings, enforcing the linear dependence of decisions on uncertain parameters. Since their first introduction in stochastic programming \cite{holt1955linear}, affine decision rules have attracted broad attention in robust optimization \cite{el2021optimality} and control theory \cite{bemporad2003min,kerrigan2003robust,skaf2010design,hadjiyiannis2011efficient} due to their superior tractability. They have also been applied in various domains, including portfolio management \cite{calafiore2008multi} and power system operation \cite{warrington2012robust,lorca2016multistage,cho2022affine}. 
Specifically, our decision rule is defined as an affine function of the current demand at any fixed time point, with coefficients invariant over time and constant terms that are continuous and piecewise affine in time. 
Exploiting the time-invariance of the coefficients and the finite number of breakpoints of the constant terms, we reformulate our model as a linear program (LP) that can be solved using the cutting-plane method \cite{kelley1960cutting}.

To the best of our knowledge, the proposed model is the first continuous-time generation scheduling model to ensure power supply-demand balance in a non-anticipative manner across an infinite number of demand trajectories. It can be generally applied to any future time horizon, such as 24 hours. However, integrating our model into current electricity market frameworks is beyond the scope of this paper.

The rest of this paper is organized as follows: 
Section \ref{sec:formulation} formulates the general robust continuous-time generation scheduling problem, introducing the decision rule framework.
Sections \ref{sec:proposed} and \ref{sec:solutionmethod} present our model and explain the solution method, respectively. 
Section \ref{sec:numericalexamples} provides numerical examples. 
Section \ref{sec:conclusions} offers 
concluding remarks.

\noindent
{\bf Notation}. We denote the sets of all, non-negative, and non-positive real numbers by $\mathbb R$, $\mathbb R_+$, and $\mathbb R_-$, respectively. For any positive integer $n$, we denote the $n$-dimensional vector of ones by ${1}_n$ and let ${\mathbb N}\left[n\right]:=\left\{1,\ldots,n\right\}$. 
Furthermore, $\circ$ represents the entrywise product operator for two vectors of the same dimension.
The vertex set of a compact convex polytope is denoted by ${\mathcal V}\left(\cdot\right)$.

\section{Problem Formulation}\label{sec:formulation}
We consider a transmission system of $G$ fossil-fuel-based generators, $D$ loads, and $L$ transmission lines for a planning horizon $T\in\left(0,\infty\right)$. 
The generators are assumed to remain in operation (i.e., turned on) throughout the planning horizon. 
For any positive integer $n$, we denote by ${\mathcal C}_n$ the set of all continuous functions from $\mathcal T:=\left[0,T\right]$ to $\mathbb R^n$. 
We also denote the demand trajectory of load $d\in{\mathbb N}\left[D\right]$ by $\xi_d:\mathcal T\to{\mathbb R}$, which is initially unknown. 
Moreover, we represent a set of possible demand trajectories of load $d$ as $\Xi_d$, which is assumed to be given and non-singleton. 
We define $\xi:{\mathcal T}\to{\mathbb R}^D$ as $\xi\left(t\right):=\left(\xi_1\left(t\right),\ldots,\xi_D\left(t\right)\right)$ and let $\Xi:=\Xi_1\times\cdots\times\Xi_D$. 
We assume that $\xi$ is revealed over time, i.e., $\Pi_t\xi:\left[0,t\right]\to{\mathbb R}^D$ is available at time $t\in{\mathcal T}$, where $\Pi_t$ is a truncation operator such that $\Pi_t\xi\left(\tau\right)=\xi\left(\tau\right)$ for any $\tau\in\left[0,t\right]$.

To ensure that power supply and demand are balanced at any time $t\in{\mathcal T}$ for any demand trajectory $\xi\in\Xi$ while satisfying the non-anticipativity constraint, we adopt a decision rule approach. Specifically, our goal is to determine a function $x:\Xi\times{\mathcal T}\to{\mathbb R}^G_+$ such that $x\left(\xi,t\right)$ denotes the vector of power outputs from the generators at time $t$ given that the demand follows $\xi$. To this end, we first consider the following general robust generation scheduling problem: 
\begin{subequations}\label{eq:prob_ori}
\begin{align}
&\inf_{x:\Xi\times{\mathcal T}\to{\mathbb R}^G_+} f_\Xi\left(x\right)\label{eq:prob_ori_obj}\\
&\text{s.t.}\quad x\left(\xi,\cdot\right)\in{\mathcal C}_{G}\quad\forall\xi\in\Xi,\label{eq:y_continuity}\\
\begin{split}
&x\left(\xi,t\right)=x\left(\xi',t\right)\ \forall(\xi,\xi',t)\in{\Xi}^2\times{\mathcal T}: 
\Pi_t\xi=\Pi_t\xi^\prime,
\end{split}\label{eq:y_na}\\
&\underline{X} \le x\left(\xi,t\right)\le \overline{X}\quad
\forall \left(\xi,t\right)\in\Xi\times{\mathcal T},\label{eq:prob_ori_con_cap}\\
&-\overline{F} \le F^{\rm g} x\left(\xi,t\right) + F^{\rm d}
\xi\left(t\right)
\le
\overline{F}
\ \ 
\forall \left(\xi,t\right)\in\Xi\times{\mathcal T},\label{eq:prob_ori_con_line}\\
\begin{split}
&\underline{R} \le \frac{x\left(\xi,{t_2}\right) - x\left(\xi,{t_1}\right)}{t_2-t_1}
\le \overline{R} \ \ \forall\left(\xi,t_1,t_2\right)\in 
\Xi\times{\mathcal T}^{\rm o}, 
\end{split}\label{eq:prob_ori_con_ramp}\\
&{1}^\top_G x\left(\xi,t\right) = {1}^\top_D\xi\left(t\right) \quad
\forall \left(\xi,t\right)\in\Xi\times{\mathcal T},\label{eq:prob_ori_con_bal}
\end{align}
\end{subequations}
which is based on a DC power flow representation \cite{stott2009dc}.

Assuming that the generation cost functions are linear in power output, we define the objective function in (\ref{eq:prob_ori_obj}) as 
\[
f_\Xi\left(x\right):=
\sup_{\xi\in\Xi} 
\int_{\mathcal T} C^\top x\left(\xi,t\right)dt
\]
where $C\in{\mathbb R}^G_+$ denotes the cost coefficient vector of the generators. Constraints (\ref{eq:y_continuity}) and (\ref{eq:y_na}) indicate that any power output trajectory must be continuous and non-anticipative, respectively. In (\ref{eq:prob_ori_con_cap}), the entries of $\overline{X},\underline{X}\in{\mathbb R}^{G}_+$ correspond to the upper and lower limits of power outputs, respectively. Constraint (\ref{eq:prob_ori_con_line}) limits the power flow in each transmission line to its respective capacity, where $\overline{F} \in \mathbb{R}^L_+$ contains the transmission capacities. The entries of $F^{\rm g}\in{\mathbb R}^{L\times G}$ and $F^{\rm d}\in{\mathbb R}^{L\times D}$ are the power transfer distribution factors of the buses to which the generators and the loads are connected, respectively. Constraint (\ref{eq:prob_ori_con_ramp}) imposes the ramp rate limits on power outputs, where 
\[
{\mathcal T}^{\rm o}:=\left\{\left(t_1,t_2\right)\in{\mathcal T}^2: t_1< t_2\right\}
\]
includes all the strictly ordered pairs of time points. 
Moreover, $\overline{R}\in{\mathbb R}^G_+$ and $\underline{R}\in{\mathbb R}^G_-$ denote the vectors of 
ramp-up and ramp-down limits of the generators, respectively. 
The systemwide supply-demand balance condition is expressed by (\ref{eq:prob_ori_con_bal}).

The most distinctive feature of (\ref{eq:prob_ori}) is the non-anticipativity constraint (\ref{eq:y_na}), which the existing models using Bernstein polynomials, e.g., \cite{le2022aggregation,khatami2019flexibility,khatami2019stochastic,hreinsson2019continuous,nikoobakht2020continuous}, lack by construction. 
However, (\ref{eq:prob_ori}) is insufficiently specific to solve directly, requiring to determine infinitely many functions from 
$\left\{\Pi_t\xi: \xi\in\Xi\right\}$ 
to ${\mathbb R}^{G}_+$ for all $t\in{\mathcal T}$. 
For tractability, we address (\ref{eq:prob_ori}) suboptimally.

\section{Proposed Approach}
\label{sec:proposed}
In this section, we explain our approach for approximately solving (\ref{eq:prob_ori}) and formulate our model. Our approach is based on the following assumption on the rate of change in demand:  
\begin{assumption}\label{ass:demand_ramp} 
For each $d\in{\mathbb N}\left[D\right]$, any demand trajectory $\xi_d\in\Xi_d$ is Lipschitz continuous, i.e., 
\[
\begin{aligned}
&\overline{r}_d:= 
\sup\left\{
R\left(\xi_d,t_1,t_2\right): 
\xi_d\in\Xi_d,\left(t_1,t_2\right)\in{\mathcal T}^{\rm o}
\right\}<\infty,\\
&\underline{r}_d:= 
\inf\left\{R\left(\xi_d,t_1,t_2\right): \xi_d\in\Xi_d,\left(t_1,t_2\right)\in{\mathcal T}^{\rm o}
\right\}>-\infty
\end{aligned}
\]
where
\[
R\left(\xi_d, t_1, t_2\right):=\frac{\xi_d\left(t_2\right)-\xi_d\left(t_1\right)}{t_2-t_1}
\]
denotes the average rate of change in the demand of load $d$ over the time interval $[t_1,t_2]$.
\end{assumption}

Our model is a surrogate of (\ref{eq:prob_ori}), formulated in consideration of tractability by replacing $\Xi$ with a superset 
$\Omega\supseteq\Xi$ and confining $x$ to a certain class $\hat{x}:{\Xi}\times{\mathcal T}\to{\mathbb R}^G_+$.

\subsection{Approximation of Demand Trajectory Sets}\label{sec:howtoomega}
In this subsection, we describe how $\Omega$ is constructed. 
We approximate $\Xi$ using only the supremum and infimum of demand for each of finitely many time intervals that subdivide $\mathcal T$, as well as the supremum and infimum of temporal changes in demand under Assumption \ref{ass:demand_ramp}. 
With this aim, we consider $N\ge2$ equally spaced time points $T_1$, \ldots, $T_N$ such that 
\[
0=T_1<T_2<\ldots<T_N=T.
\] 
The number $N$ of time points affects both conservatism of $\Omega$ and the computational burden of our model, which should be set taking into account their trade-off. 
Then, the superset $\Omega$ of $\Xi$ is defined as $\Omega:=\Omega_1\times\cdots\times\Omega_D
\subset{\mathcal C}_D$, where 
\[
\begin{aligned}
&{\Omega}_d:= 
\{
\xi_d
\in{\mathcal C}_1:
\underline{\omega}^\prime_d\left(t\right)\le\xi_d\left(t\right)\le\overline{\omega}^\prime_d\left(t\right)\quad\forall t\in{\mathcal T},\\ 
&\quad\underline{r}_d \le 
R\left(\xi_d,t_1,t_2\right)
\le\overline{r}_d
\quad\forall\left(t_1,t_2\right)\in{\mathcal T}^{\rm o}\}.
\end{aligned}
\]
Here, we let 
\[
\begin{aligned}
&\overline{\omega}^\prime_d\left(t\right):=
\sup\left\{\xi_d\left(\tau\right): 
\xi_d\in\Xi_d,\tau\in\left[T_{i\left(t\right)},T_{i\left(t\right)+1}\right]\right\}
\\ 
&\underline{\omega}^\prime_d\left(t\right):=
\inf\left\{\xi_d\left(\tau\right): 
\xi_d\in\Xi_d,\tau\in\left[T_{i\left(t\right)},T_{i\left(t\right)+1}\right]\right\}
\end{aligned}
\]
with $i:{\mathcal T}\to{\mathbb N}\left[N-1\right]$ defined as 
\[
i\left(t\right):=
\begin{cases}
\begin{aligned}
&j, &&\text{if}\quad 
t\in\left[T_j,T_{j+1}\right),\\
&N-1 && \text{if}\quad t=T.
\end{aligned} 
\end{cases} 
\] 
The index function $i$ returns the time interval to which each time point belongs, out of the initially given $\left(N-1\right)$ time intervals. Thus, $\overline{\omega}^\prime_d$ and $\underline{\omega}^\prime_d$ are the tightest piecewise constant functions that capture $\Xi_d$, each defined with $\left(N-1\right)$ segments corresponding to the $\left(N-1\right)$ time intervals $\left[T_1,T_2\right),\ldots,\left[T_{N-2},T_{N-1}\right),\left[T_{N-1},T_N\right]$.
Note that we further constrain the ramp rates of load $d$'s demand using $\overline{r}_d$ and $\underline{r}_d$ to define $\Omega_d$, in accordance with Assumption \ref{ass:demand_ramp}.

Our first step in addressing (\ref{eq:prob_ori}) is to replace $\Xi$ by $\Omega$, which, along with $\hat{x}$, reduces our model to a finite-dimensional program. This is attributed to the fact that $\Omega_d$  for any $d\in{\mathbb N}\left[D\right]$ can be represented using two piecewise affine functions $\overline{\omega}_d,\underline{\omega}_d\in{\mathcal C}_1$, 
which are defined in what follows.

A continuous piecewise affine function over a closed interval can be uniquely identified by its breakpoints and the corresponding function values. Thus, to define $\overline\omega_d$ and $\underline\omega_d$, we present their function values at breakpoints.

We first define $\overline{\omega}_d$ for any fixed $d\in{\mathbb N}\left[D\right]$. The breakpoint set of $\overline\omega_d$ is denoted by
\[
\overline{\mathcal T}^{\rm b}_d:=\left\{T_1,\ldots,T_N\right\}\cup
\{\overline{T}^{\rm b}_{d,1},\ldots,\overline{T}^{\rm b}_{d,N-2}\}
\]
where
\[
\overline{T}^{\rm b}_{d,j}:= 
\begin{cases}
\begin{aligned}
&T_{j+1} + 
\overline{\omega}^{\prime+}_{d,j}/\overline{r}_d
\quad\text{if}\quad
\overline{\omega}^{\prime+}_{d,j}
\ge0,\\
&T_{j+1} - \overline{\omega}^{\prime+}_{d,j}/\underline{r}_d
\quad \text{otherwise}
\end{aligned}
\end{cases}
\]
with 
\[
\overline{\omega}^{\prime+}_{d,j} :=\overline\omega^\prime_d\left(T_{j+1}\right) - \overline\omega^\prime_d\left(T_j\right).
\]
Here, we follow the convention $0/0=0$. Note that $\overline{T}^{\rm b}_{d,j}\in\left(T_{j},T_{j+2}\right)$ for any $j\in{\mathbb N}\left[N-2\right]$. 
Meanwhile, the function values of $\overline\omega_d$ for each $t\in\overline{\mathcal T}^{\rm b}_d$ are 
defined as follows: 
\[
\begin{aligned}
&\overline{\omega}_d\left(T_1\right):= \overline{\omega}^\prime_d\left(T_1\right),\\
&\overline{\omega}_d\left(T_{j+1}\right):=
\min\left\{
\overline{\omega}^\prime_d\left(T_{j}\right),\overline{\omega}^\prime_d\left(T_{j+1}\right)
\right\}&&\forall j\in{\mathbb N}\left[N-1\right],\\
&\overline{\omega}_d\,\big(\overline{T}^{\rm b}_{j}\big) := 
\max\left\{
\overline{\omega}^\prime_d\left(T_{j}\right),\overline{\omega}^\prime_d\left(T_{j+1}\right)
\right\}&&\forall j\in{\mathbb N}\left[N-2\right],
\end{aligned}
\]
whereby $\overline\omega_d$ is fully identified.
 
Intuitively, $\overline{\omega}_d$ is obtained by taking the 
highest possible ramp rate at each time point given the initial value $\overline{\omega}_d\left(T_1\right)$ and the piecewise constant functions $\overline{\omega}^\prime_d$ and $\underline{\omega}^\prime_d$, while ensuring that 
$\overline{\omega}_d\in\Omega_d$. As such, $\overline{\omega}_d$ is the pointwise largest function in $\Omega_d$, i.e., 
\[
\overline{\omega}_d\left(t\right)
=
\sup\left\{
\xi_d\left(t\right): \xi_d\in\Omega_d\right\}\quad\forall t\in{\mathcal T}.
\]

\begin{figure}[t!]
    \centering
       \includegraphics[width=8.5cm,height=2.5cm]
       {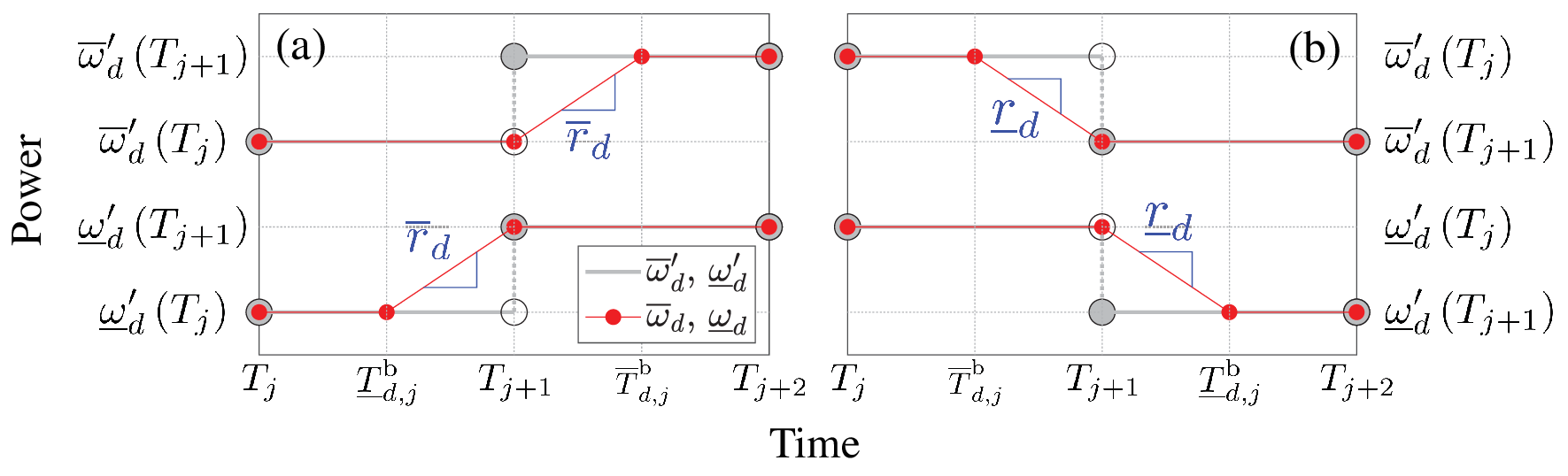}
    \caption{Obtaining $\left(\overline{\omega}_d,\underline{\omega}_d\right)$.
    (a) $\overline{\omega}^{\prime+}_{d,j},\underline{\omega}^{\prime+}_{d,j}\ge0$. 
    (b) $\overline{\omega}^{\prime+}_{d,j},\underline{\omega}^{\prime+}_{d,j}\le0$.} 
\label{fig:1}
\end{figure}

We define $\underline{\omega}_d$ in a similar manner. Specifically, we let
\[
\underline{\mathcal T}^{\rm b}_d:=\left\{T_1,\ldots,T_N\right\}\cup\{\underline{T}^{\rm b}_{d,1},\ldots,\underline{T}^{\rm b}_{d,N-2}\}
\]
denote the breakpoint set of $\underline\omega_d$, where 
\[
\underline{T}^{\rm b}_{d,j} := 
\begin{cases}
\begin{aligned}
&T_{j+1} - 
\underline{\omega}^{\prime+}_{d,j}/\overline{r}_d
\quad\text{if}\quad
\underline{\omega}^{\prime+}_{d,j}\ge0,\\
&T_{j+1} + \underline{\omega}^{\prime+}_{d,j}/\underline{r}_d
\quad\text{otherwise}
\end{aligned}
\end{cases}
\]
with 
\[
\underline{\omega}^{\prime+}_{d,j}:=\underline\omega^\prime_d\left(T_{j+1}\right) - \underline\omega^\prime_d\left(T_j\right).
\] 
Moreover, the function values of $\underline{\omega}_d$ for each $t\in\underline{\mathcal T}^{\rm b}_d$ are 
defined as follows: 
\[
\begin{aligned}
&\underline{\omega}_d\left(T_1\right):= \underline{\omega}^\prime_d\left(T_1\right),\\
&\underline{\omega}_d\left(T_{j+1}\right):=
\max\left\{
\underline{\omega}^\prime_d\left(T_{j}\right),\underline{\omega}^\prime_d\left(T_{j+1}\right)
\right\}&&\forall j\in{\mathbb N}\left[N-1\right],\\
&\underline{\omega}_d\,\big(\underline{T}^{\rm b}_{j}\big):= 
\min\left\{
\underline{\omega}^\prime_d\left(T_{j}\right),\underline{\omega}^\prime_d\left(T_{j+1}\right)
\right\}&&\forall j\in{\mathbb N}\left[N-2\right]. 
\end{aligned}
\]
Given its initial value $\underline{\omega}_d\left(T_1\right)$, 
$\underline{\omega}_d$ follows the lowest possible ramp rate at each time point while ensuring that $\underline{\omega}_d\in\Omega_d$. Thus, $\underline{\omega}_d$ is the pointwise smallest function in $\Omega_d$, i.e.,
\[
\underline{\omega}_d\left(t\right)
=
\inf\left\{
\xi_d\left(t\right): \xi_d\in\Omega_d\right\}\quad\forall t\in{\mathcal T}.
\]
Fig. \ref{fig:1}
depicts how to obtain $\left(\overline{\omega}_d,\underline{\omega}_d\right)$ using $\left(\overline{\omega}^\prime_d,\underline{\omega}^\prime_d,\overline{r}_d,\underline{r}_d\right)$.

Notably, the continuous piecewise affine functions $\overline{\omega}_d$ and $\underline{\omega}_d$ offer the tightest interval representation of $\Omega_d$ in the sense that 
\[ 
\left\{\xi_d\left(t\right): \xi_d\in\Omega_d\right\} = \left[\underline{\omega}_d\left(t\right),\overline{\omega}_d\left(t\right)\right]
\]
for any $t\in{\mathcal T}$. 
We define $\overline{\omega},\underline{\omega}:{\mathcal T}\to{\mathbb R}^D$ 
as 
\[
\overline{\omega}\left(t\right):=\left(\overline{\omega}_1(t),\ldots,\overline{\omega}_D(t)\right),\quad\underline{\omega}\left(t\right):=\left(\underline{\omega}_1(t),\ldots,\underline{\omega}_D(t)\right).
\]

\subsection{Affine Decision Rule}
In this subsection, we present the proposed decision rule $\hat{x}$. 
Let the set of all breakpoints of $\overline\omega$ and $\underline\omega$ be denoted 
by 
\[
{\mathcal T}^{\rm b}:=
\cup_{d\in{\mathbb N}\left[D\right]}
\,\big(\overline{\mathcal T}^{\rm b}_d\cup\underline{\mathcal T}^{\rm b}_d\big) = \left\{T^{\rm b}_1,T^{\rm b}_2,\ldots,T^{\rm b}_M\right\}
\]
where 
\[
0=T^{\rm b}_1<T^{\rm b}_2<\ldots<T^{\rm b}_M=T
\] 
and $M=\lvert{\mathcal T}^{\rm b}\rvert\le 2D\left(N-2\right) + N$. 
We define 
$k: {\mathcal T}
\to{\mathbb N}\left[
M-1
\right]$ 
as 
\[
k\left(t\right):=
\begin{cases}
\begin{aligned}
& j &&\text{if}\quad 
t\in\left[T^{\rm b}_j,T^{\rm b}_{j+1}\right),\\
&M-1 &&\text{if}\quad t=T.
\end{aligned}
\end{cases}
\]
Similar to $i$, $k$ corresponds to the index function for the $\left(M-1\right)$ time intervals formed by the $M$ time points in ${\mathcal T}^{\rm b}$. We also define 
$\gamma:{\mathcal T}\to\left[0,1\right]$ as 
\[
\gamma\left(t\right) := 
\begin{cases}
\begin{aligned}
&(t-T^{\rm b}_{k\left(t\right)})/
\Delta^{\rm b}_{k\left(t\right)}, &&\text{if}\quad t\in\left[0,T\right), \\
&1, &&
\text{if}\quad t=T,
\end{aligned}
\end{cases}
\]
where 
\[
\Delta^{\rm b}_j:=T^{\rm b}_{j+1} - T^{\rm b}_{j}\quad\forall j\in{\mathbb N}\left[M-1\right].
\]
Note that $\gamma\left(t\right)$ represents the position of time point $t$ within the time interval $[T^{\rm b}_{k(t)},T^{\rm b}_{k(t)+1}]$ to which it belongs. As a result, we can express any $t\in{\mathcal T}$ using $k$ and $\gamma$ 
as 
\[
t=T^{\rm b}_{k\left(t\right)} + \gamma\left(t\right)\Delta^{\rm b}_{k\left(t\right)}.
\] 
Then, our decision rule $\hat{x}:\Xi\times{\mathcal T}\to{\mathbb R}^G_+$ is defined as 
\[
\hat{x} \left(\xi,t\right):= \alpha\xi\left(t\right) + 
\left(1-\gamma\left(t\right)\right) \beta_{k\left(t\right)} + \gamma\left(t\right)\beta_{k\left(t\right)+1}
\]
where $\alpha\in{\mathbb R}^{G\times D}$ and $\beta_1,\ldots,\beta_M\in{\mathbb R}^{G}$ parameterize $\hat{x}$. 
For brevity, we define $\beta\in{\mathcal C}_G$ as 
\[
\beta\left(t\right):= \left(1-\gamma\left(t\right)\right) \beta_{k\left(t\right)} + \gamma\left(t\right)\beta_{k\left(t\right)+1},
\] 
which is continuous and piecewise affine in time, and $\beta_{1:M}:=\left(\beta_1,\ldots,\beta_M\right)$. 
When needed, we use the symbol $\hat{x}_{\alpha,\beta_{1:M}}$ instead of $\hat{x}$ to explicitly express its dependency on $\left(\alpha,\beta_{1:M}\right)$.

Our decision rule is analogous to affine decision rules widely studied in discrete-time settings, e.g., \cite{lorca2016multistage}, as it remains affine in the uncertain parameters at any fixed time point. Specifically, $\hat{x}$ can be viewed as the linear interpolation of affine decision rules defined for each time point in ${\mathcal T}^{\rm b}$ with a common coefficient matrix. We present our generation scheduling model in the following subsection.

\subsection{Proposed Model}
By replacing $\Xi$ and $x$ in (\ref{eq:prob_ori}) by $\Omega$ and $\hat{x}$, respectively, we formulate our model to determine $\left(\alpha,\beta_{1:M}\right)$ as 
\begin{subequations}\label{eq:prob_our}
\begin{align}
\inf_{\alpha,\beta_{1:M}} &f_\Omega\left(\hat{x}_{\alpha,\beta_{1:M}}\right)
\label{eq:prob_our_obj}\\
\text{s.t.}
\ \ &A\hat{x}_{\alpha,\beta_{1:M}}\left(\xi,t\right)
+B\xi\left(t\right) \le a\ 
\forall \left(\xi,t\right)\in\Omega\times{\mathcal T},\label{eq:prob_our_con_ineq}\\
\begin{split}
&\underline{R} \le \frac{\hat{x}_{\alpha,\beta_{1:M}}\left(\xi,{t_2}\right) - \hat{x}_{\alpha,\beta_{1:M}}\left(\xi,{t_1}\right)}{t_2-t_1}
\le \overline{R} \\
&\quad \forall\left(\xi,t_1,t_2\right)\in 
\Omega\times{\mathcal T}^{\rm o}, 
\end{split}\label{eq:prob_our_con_ramp}\\
&{1}^\top_G \hat{x}_{\alpha,\beta_{1:M}}\left(\xi,t\right) = {1}^\top_D\xi\left(t\right) \quad
\forall \left(\xi,t\right)\in\Omega\times{\mathcal T},\label{eq:prob_our_con_bal}\\
&\alpha\in{\mathbb R}^{G\times D},\beta_{1:M}\in{\mathbb R}^{GM}\nonumber
\end{align}
\end{subequations}
where $A\in{\mathbb R}^{2(G+L)\times G}$, $B\in{\mathbb R}^{2(G+L)\times D}$, and $a\in{\mathbb R}^{2(G+L)}$ are defined using $\overline{X}$, $\underline{X}$, $F^{\rm g}$, $F^{\rm d}$, and $\overline{F}$ such that (\ref{eq:prob_our_con_ineq}) is collectively equivalent to (\ref{eq:prob_ori_con_cap}) and (\ref{eq:prob_ori_con_line}) for $\left(\Xi,x\right)=\left(\Omega,\hat{x}\right)$. Moreover, we omit the continuity constraint (\ref{eq:y_continuity}) and the non-anticipativity constraint (\ref{eq:y_na}) for $\left(\Xi,x\right)=\left(\Omega,\hat{x}\right)$ as our decision rule satisfies both constraints by construction. We assume that (\ref{eq:prob_our}) is feasible. 
Notably, $\hat{x}$ specified by any feasible point of (\ref{eq:prob_our}) meets all the constraints of (\ref{eq:prob_ori}). Further, the optimal value of (\ref{eq:prob_our}) is an upper bound on that of (\ref{eq:prob_ori}). We discuss a solution method for (\ref{eq:prob_our}) in the next section.

\section{Solution Method} \label{sec:solutionmethod} 
In this section, we show that (\ref{eq:prob_our}) can be solved using the cutting plane method in finitely many iterations. First, we present an LP reformulation of (\ref{eq:prob_our}). Subsequently, we describe the cutting-plane method for solving this LP.

\subsection{LP Reformulation}
To solve (\ref{eq:prob_our}), we first reformulate it as an LP by expressing the objective function in (\ref{eq:prob_our_obj}) as a convex, piecewise affine function of $\left(\alpha,\beta_{1:M}\right)$ and deriving a finite set of linear inequalities for $\left(\alpha,\beta_{1:M}\right)$ corresponding to (\ref{eq:prob_our_con_ineq})--(\ref{eq:prob_our_con_bal}). 
These equivalents of (\ref{eq:prob_our_obj})--(\ref{eq:prob_our_con_bal}) are provided in Propositions \ref{prop:cost}--\ref{prop:bal}, respectively.

\begin{proposition}\label{prop:cost}
The objective function in 
(\ref{eq:prob_our_obj}) can be rewritten as 
follows: 
\[
\begin{aligned}
&f_\Omega\left(\hat{x}_{\alpha,\beta_{1:M}}\right)
= 
\frac{1}{2}
\sum_{j=1}^{M-1}\Delta^{\rm b}_j C^\top(\beta_{j+1}+\beta_j)\\
&\quad+\frac{1}{2}\sum_{d=1}^D\sup\left\{
\sum_{j=1}^{M-1} \Delta^{\rm b}_j 
\overline{\omega}^{\rm s}_{d,j}
C^\top\alpha_d, 
\sum_{j=1}^{M-1} \Delta^{\rm b}_j \underline{\omega}^{\rm s}_{d,j}
C^\top \alpha_{d}\right\}
\end{aligned}
\]
where $\alpha_{d}\in{\mathbb R}^G$ denotes the $d$th column of $\alpha$ and
\[
\begin{aligned}
&\overline{\omega}^{\rm s}_{d,j} := \overline\omega_d\left(T^{\rm b}_{j+1}\right) + \overline\omega_d\left(T^{\rm b}_{j}\right), \\
&\underline{\omega}^{\rm s}_{d,j}:= \underline\omega_d\left(T^{\rm b}_{j+1}\right) + \underline\omega_d\left(T^{\rm b}_{j}\right).
\end{aligned}
\]
\end{proposition}
The proof of Proposition \ref{prop:cost} is deferred to Appendix \ref{app:prop:cost}.

Due to the time-invariance of $\alpha$ and the linearity of the generation cost of each generator,
the total generation cost for any fixed $(\alpha,\beta_{1:M})$ can be represented as an affine function of the total demand of each load, i.e., the integral of $\xi$, whose constant term includes the integral of $\beta$. This implies that the pointwise maximum or minimum trajectory, i.e., $\overline{\omega}_d$ or $\underline{\omega}_d$, corresponds to the worst-case scenario of $\xi_d$ for any $d\in{\mathbb N}\left[D\right]$ in terms of total generation cost. Note that all entries of $\overline\omega$, $\underline\omega$, and $\beta$ are piecewise affine functions. Moreover, the integral of a piecewise affine function over an interval can be computed by summing the areas of the trapezoids associated with each segment. Proposition \ref{prop:cost} follows as a direct consequence of this observation.

\begin{proposition}\label{prop:ineq}
Let $\left(\sigma_{1},\ldots,\sigma_{2^D}\right)$ denote any permutation of $\left\{0,1\right\}^D$. Moreover, for each $j\in{\mathbb N}\left[M\right]$, let 
\[
H_{j} :=\overline{\omega}\left(T^{\rm b}_j\right)-\underline{\omega}\left(T^{\rm b}_j\right).
\] 
Constraint (\ref{eq:prob_our_con_ineq}) is necessary and sufficient for 
\begin{equation}\label{eq:alphabeta_both_cap_lines}
\begin{aligned}
&\left(A\alpha + B \right) \left(\underline{\omega}\left(T^{\rm b}_j\right) + \sigma_l \circ
H_{j}
\right) + A\beta_j \le a\\
&\quad\forall \left(j,l\right)\in{\mathbb N}\left[M\right]\times{\mathbb N}\left[2^D\right].
\end{aligned}
\end{equation}
\end{proposition}
The proof of Proposition \ref{prop:ineq} is provided in Appendix \ref{app:prop:ineq}.

Note that $(\underline{\omega}\left(T^{\rm b}_j\right) + \sigma_l \circ H_{j})$ in (\ref{eq:alphabeta_both_cap_lines}) represents a vertex of the multidimensional interval $[\underline{\omega}(T^{\rm b}_j),\overline{\omega}(T^{\rm b}_j)]$. Thus, (\ref{eq:alphabeta_both_cap_lines}) ensures that our decision rule satisfies the inequality in (\ref{eq:prob_our_con_ineq}) at each time point in ${\mathcal T}^{\rm b}$ for all extremal demand scenarios. This, in turn, implies that any inequality in (\ref{eq:prob_our_con_ineq}) can be obtained as a convex combination of the inequalities in (\ref{eq:alphabeta_both_cap_lines}) because $\overline\omega$, $\underline\omega$, and $\beta$ are piecewise affine functions of time with the same breakpoint set $\mathcal T^{\rm b}$.

\begin{proposition}\label{prop:ramp}
For each $\left(d,j\right)\in{\mathbb N}\left[D\right]\times{\mathbb N}\left[M-1\right]$, if 
\[
\overline{\omega}_d (T^{\rm b}_j)=\underline{\omega}_{d} (T^{\rm b}_j),\quad\overline{\omega}_d (T^{\rm b}_{j+1})=\underline{\omega}_{d} (T^{\rm b}_{j+1}),
\]
implying that load $d$ has only a single demand trajectory over the time interval $[T^{\rm b}_{j},T^{\rm b}_{j+1}]$, then let 
\[
\begin{aligned}
&\overline{r}^\prime_{d,j}:=\left(\overline{\omega}_d(T^{\rm b}_{j+1})-\overline{\omega}_d(T^{\rm b}_{j})\right)/\Delta^{\rm b}_j,\\ 
&\underline{r}^\prime_{d,j}:=\left(\overline{\omega}_d(T^{\rm b}_{j+1})-\overline{\omega}_d(T^{\rm b}_{j})\right)/\Delta^{\rm b}_j.
\end{aligned}
\]
Otherwise, let $\overline{r}^\prime_{d,j}:=\overline{r}_d$ and $\underline{r}^\prime_{d,j}:=\underline{r}_d$. 
Constraint (\ref{eq:prob_our_con_ramp}) is sufficient and necessary for
\begin{equation}\label{eq:alphabeta_ramp}
\begin{aligned}
&\underline{R}\le\alpha
\left(
\underline{r}^\prime_j + \sigma_l\circ\left(\overline{r}^\prime_j-\underline{r}^\prime_j\right)\right) + g_j \le \overline{R}\\
&\quad \forall \left(j,l\right)\in{\mathbb N}\left[M-1\right]\times{\mathbb N}\left[2^D\right]
\end{aligned}
\end{equation}
where $\underline{r}^\prime_j:=(\underline{r}^\prime_{j,1},\ldots,\underline{r}^\prime_{j,D})\in{\mathbb R}^D$, 
$\overline{r}^\prime_j:=(\overline{r}^\prime_{j,1},\ldots,\overline{r}^\prime_{j,D})\in{\mathbb R}^D$, 
and 
\[
g_j:=(\beta_{j+1} - \beta_j)/\Delta^{\rm b}_j\in{\mathbb R}^G.
\] 
\end{proposition}
The proof of Proposition \ref{prop:ramp} is given in Appendix \ref{app:prop:ramp}.

In (\ref{eq:alphabeta_ramp}), $\overline{r}^\prime_j$ and $\underline{r}^\prime_j$ denote the maximum and minimum ramp rates of each load's demand over the time interval $[T^{\rm b}_j,T^{\rm b}_{j+1}]$, respectively.  
Moreover, $(\underline{r}^\prime_j + \sigma_l \circ (\overline{r}^\prime_j - \underline{r}^\prime_j))$ represents a vertex of the multidimensional interval $[\underline{r}^\prime_j,\overline{r}^\prime_j]$. Thus, (\ref{eq:alphabeta_ramp}) constrains $\alpha$ and the rate of change in $\beta$ for each time interval defined by the breakpoints of $\beta$, ensuring that the ramp rate limits of the generators are satisfied for any sudden fluctuation of demand. Since the ramp rate of a generator's power output trajectory obtained using our decision rule is determined by the ramp rates of $\xi$ and $\beta$, the equivalence of (\ref{eq:prob_our_con_ramp}) and (\ref{eq:alphabeta_ramp}) is intuitive.

\begin{proposition}\label{prop:bal}
Constraint (\ref{eq:prob_our_con_bal}) is sufficient and necessary for 
\begin{equation}\label{eq:alphabeta_eq}
1^\top_G\alpha = {1}^\top_D,\quad
{1}^\top_G\beta_j = 0\quad\forall j\in{\mathbb N}\left[M\right]. 
\end{equation}
\end{proposition}
The proof of Proposition \ref{prop:bal} is presented in Appendix \ref{app:prop:bal}.

Given that $\Xi_d$ is non-singleton for any $d\in{\mathbb N}\left[D\right]$ by assumption, Proposition \ref{prop:bal} is essential. It is worth mentioning that the same result generally holds for affine decision rules in discrete-time generation scheduling models (e.g., \cite{lorca2016multistage}).

Using Propositions \ref{prop:cost}--\ref{prop:bal}, (\ref{eq:prob_our}) is recast as the LP 
\begin{equation}\label{eq:prob_our_LP}
\begin{aligned}
\inf_{\alpha,\beta_{1:M},\eta_{1:D}}
\ &\frac{1}{2}
\sum_{j=1}^{M-1}
\Delta^{\rm b}_{j}
C^\top\left(\beta_{j}+\beta_{j+1}\right)
+
\frac{1}{2}\sum_{d=1}^D\eta_d \\
\text{s.t.}\qquad&\text{(\ref{eq:alphabeta_both_cap_lines}), (\ref{eq:alphabeta_ramp}), (\ref{eq:alphabeta_eq})},\\
&\eta_d\ge\sum_{j=1}^{M-1}\Delta^{\rm b}_{j}\overline{\omega}^{\rm s}_{d,j}C^\top\alpha_d
\quad\forall d\in{\mathbb N}\left[D\right],\\
&\eta_d\ge\sum_{j=1}^{M-1}\Delta^{\rm b}_{j}\underline{\omega}^{\rm s}_{d,j}C^\top\alpha_d
\quad\forall d\in{\mathbb N}\left[D\right],\\
&\alpha\in{\mathbb R}^{G\times D},\beta_{1:M}\in{\mathbb R}^{GM},\eta_{1:D}\in{\mathbb R}^D
\end{aligned}
\end{equation}
where $\eta_{1:D}:=(\eta_1,\ldots,\eta_D)\in{\mathbb R}^D$. 
Although written as an LP, (\ref{eq:prob_our_LP}) may be hard to tackle in the current form due to the large number of constraints in
(\ref{eq:alphabeta_both_cap_lines}) and  (\ref{eq:alphabeta_ramp}). Specifically, 
(\ref{eq:alphabeta_both_cap_lines}) and  (\ref{eq:alphabeta_ramp}) contain 
$2^DM$ and $2^{D+1}\left(M-1\right)$ linear inequalities, respectively. To efficiently solve (\ref{eq:prob_our_LP}), we employ the cutting-plane method \cite{kelley1960cutting}.

\subsection{Cutting-Plane Method}\label{sec:cutting}

In this subsection, we describe the cutting-plane method for (\ref{eq:prob_our_LP}). For initialization, we let $\zeta_j:=\underline{\omega}\left(T^{\rm b}_{j}\right) + \sigma_{l^\prime} \circ
H_j$ for each $j\in{\mathbb N}\left[M\right]$ and
$\psi_j:=\underline{r}^\prime_j + \sigma_{l^\prime}\circ\left(\overline{r}^\prime_j-\underline{r}^\prime_j\right)$ for each $j\in{\mathbb N}\left[M-1\right]$ with a randomly chosen 
$l^{\prime}\in{\mathbb N}\left[2^D\right]$. 
Further, we solve the following LP, which is a relaxation of (\ref{eq:prob_our_LP}): 
\begin{equation}\label{eq:prob_our_LP_initial}
\begin{aligned}
\inf_{\alpha,\beta_{1:M},\eta_{1:D},\eta}
\ &\eta\\
\text{s.t.}\ \qquad&\hspace{-5mm}\text{(\ref{eq:alphabeta_eq})},\\
&\hspace{-5mm}\left(A\alpha + B \right)
\zeta_j + A\beta_j \le a\quad\forall j\in{\mathbb N}\left[M\right],\\
&\hspace{-5mm}
\underline{R}
\le
\alpha
\psi_j + g_j \le \overline{R}\quad \forall j\in{\mathbb N}\left[M-1\right],\\
&\hspace{-5mm}
\eta\ge\frac{1}{2}\sum_{j=1}^{M-1} \Delta^{\rm b}_{j} C^\top\left(\beta_{j}+\beta_{j+1}\right) + \frac{1}{2}\sum_{d=1}^D\eta_d,\\
&\hspace{-5mm}
\eta_d\ge\sum_{j=1}^{M-1}\Delta^{\rm b}_{j}\underline{\omega}^{\rm s}_{d,j}C^\top\alpha_d
\quad\forall d\in{\mathbb N}\left[D\right],\\
&\hspace{-5mm}
\eta_d\ge\sum_{j=1}^{M-1}\Delta^{\rm b}_{j}\overline{\omega}^{\rm s}_{d,j}C^\top\alpha_d
\quad\forall d\in{\mathbb N}\left[D\right],\\
&\hspace{-5mm}
\alpha\in{\mathbb R}^{G\times D},\beta_{1:M}\in{\mathbb R}^{GM},\eta_{d}\in{\mathbb R},\eta\in{\mathbb R}_+
\end{aligned}
\end{equation}
where the inequalities in (\ref{eq:alphabeta_both_cap_lines}) and  (\ref{eq:alphabeta_ramp}) are removed compared to (\ref{eq:prob_our_LP}) except those for $l=l^\prime$. The new decision variable $\eta\in{\mathbb R}_+$ is added to ensure that (\ref{eq:prob_our_LP_initial}) is bounded, which does not affect the optimality as the total generation cost is known to be non-negative in any case.

After obtaining a solution $(\alpha^{(1)},\beta^{(1)}_{1:M})$ to (\ref{eq:prob_our_LP_initial}), we investigate if any inequality in (\ref{eq:alphabeta_both_cap_lines}) and (\ref{eq:alphabeta_ramp}) is violated for $(\alpha,\beta_{1:M})=(\alpha^{(1)},\beta^{(1)}_{1:M})$  with respect to $l$ for each fixed $j$. 
This can be readily achieved considering only the positivity of the entries in $A\alpha^{(1)}+B$ and $\alpha^{(1)}$, without solving an optimization problem. 
For example, to check whether the $p$th inequality in (\ref{eq:alphabeta_both_cap_lines}) is violated for some fixed $\left(p,j\right)\in{\mathbb N}\left[2\left(G+L\right)\right]\times{\mathbb N}\left[M\right]$ and for any $l\in{\mathbb N}\left[2^D\right]$, i.e., whether 
\begin{equation}\label{eq:osifs100}
\begin{aligned}
&\big(A_p\alpha^{(1)}+B_p\big)\left(
\underline{\omega}\left(T^{\rm b}_j\right) + 
\sigma_l\circ H_j
\right) + A_p\beta^{(1)}_j\le a_p\\
&\quad\forall l\in{\mathbb N}\left[2^D\right]
\end{aligned}
\end{equation}
holds, where $A_p\in{\mathbb R}^{1\times G}$,  $B_p\in{\mathbb R}^{1\times D}$, and $a_p\in{\mathbb R}$ denote the $p$th rows of $A$, $B$, and $a$,  respectively, 
it suffices to verify that 
\begin{equation}\label{eq:osifs99} 
\big(A_p\alpha^{(1)}+B_p\big)\left(
\underline{\omega}\left(T^{\rm b}_j\right) + 
\sigma^\prime\circ H_j
\right) + A_p\beta^{(1)}_j\le a_p
\end{equation} 
because we have 
\[
\max_{l\in{\mathbb N}\left[2^D\right]}\big(A_p\alpha^{(1)}+B_p\big)\left(\sigma_l\circ H_j\right)=\big(A_p\alpha^{(1)}+B_p\big)\left(\sigma^\prime\circ H_j\right).
\]
Here, as $H_j\in{\mathbb R}^D_+$, the $d$th entry of $\sigma^\prime\in\left\{0,1\right\}^D$ 
is set to one if the $d$th entry of $\left(A_p\alpha^{(1)}+B_p\right)$ is non-negative, and zero otherwise. 
If (\ref{eq:osifs99}) holds, then (\ref{eq:osifs100}) also holds. 
Otherwise, (\ref{eq:osifs99}) can be considered the most violated inequality in (\ref{eq:osifs100}). 

Subsequently, the most violated inequalities in (\ref{eq:alphabeta_both_cap_lines}) and (\ref{eq:alphabeta_ramp}), obtained in this manner, are added as constraints to (\ref{eq:prob_our_LP_initial}), which is then solved again to obtain a point $(\alpha^{(2)},\beta^{(2)}_{1:M})$. We repeat this process, meaning that at iteration $k\ge1$, we obtain $(\alpha^{(k)},\beta^{(k)}_{1:M})$ by solving (\ref{eq:prob_our_LP_initial}) with all the most violated inequalities previously found in (\ref{eq:alphabeta_both_cap_lines}) and (\ref{eq:alphabeta_ramp}) incorporated, which is an LP, until no inequality is found to be violated. Notably, the algorithm returns a solution to (\ref{eq:prob_our_LP}) in at most $2^D\left(3M-2\right)$ iterations. We omit further details of the algorithm, as it is considered standard in the literature of robust optimization; see, e.g., \cite{lorca2016multistage}.

\section{Illustrative Examples}\label{sec:numericalexamples}
In this section, we illustrate our approach with simulation results on a modified version of the six-bus test system from \cite{IIT}. The modified system consists of three generators $g=1,2,3$ at buses 1, 2 and 6, respectively, three loads $d=1, 2, 3$ at buses 3, 4 and 5, respectively, and seven transmission lines. With $T=24$ and $N=25$, we assumed that $\Xi_d$ was given for each $d$ such that $\overline\omega^\prime_d(T_{j})=\rho_d S_{j}+5$ for any $j\in{\mathbb N}\left[N-1\right]$ and $\overline\omega^\prime_d(T_N)=\overline\omega^\prime_d(T_{N-1})$, where $\rho_1=0.2$ and $\rho_2=\rho_3=0.4$. Here, $S_j$ denotes the system demand at hour $j$ reported in \cite{IIT}. 
To ensure that $\Omega$ is non-empty, we also set $\underline\omega^\prime_d(T_1)=\rho_d S_1-5$, $\underline\omega^\prime_d(T_{j})=\min \{ \overline\omega^\prime_d(T_{j-1}),  \overline\omega^\prime_d(T_{j+1}), \rho_d S_j\}-5$ for $j=2,\ldots,N-1$, and $\underline\omega^\prime_d(T_N)= \underline\omega^\prime_d(T_{N-1})$ with 
$\overline{r}_d=\rho_d\max_{j}\{S_{j+1} - S_j\}$ and $\underline{r}_d=\rho_d\min_{j}\{S_{j+1} - S_j\}$. Further, we doubled the capacities of transmission lines and halved the ramp-up and ramp-down limits of generators.

\begin{figure}[t!]
    \centering
       \includegraphics[height=3cm]
       {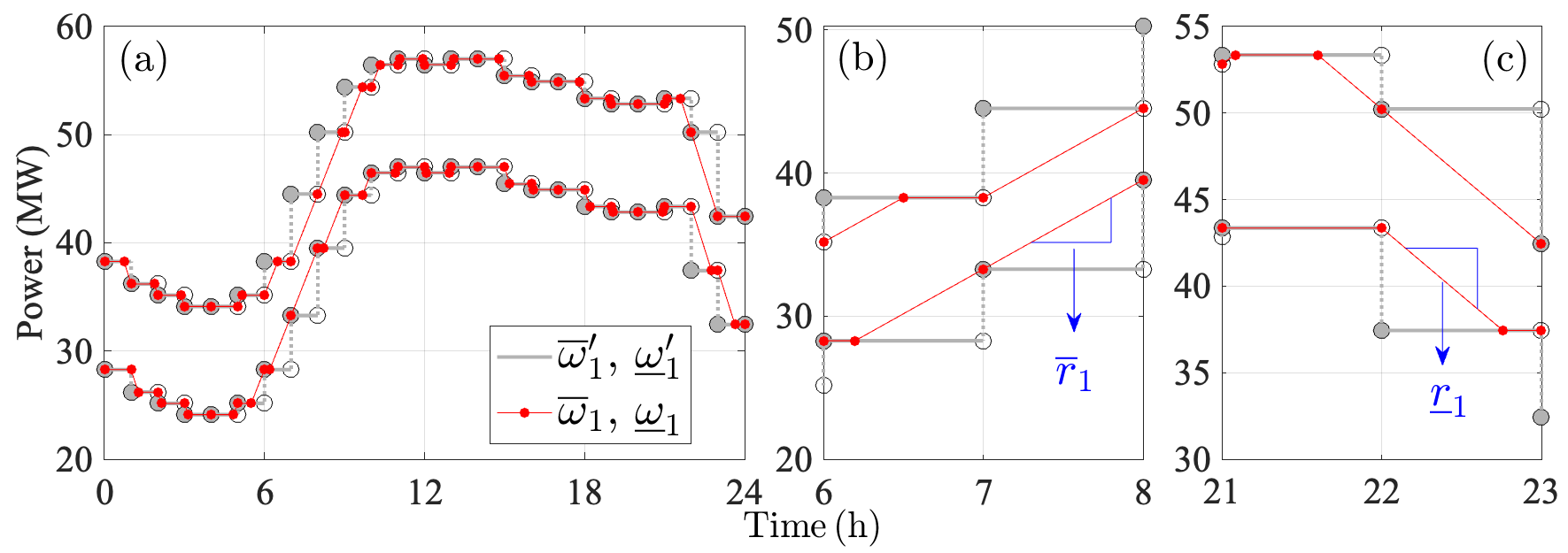}
    \caption{(a) Construction of $\Omega_1$. (b), (c) Enlarged views of (a).} 
\label{fig:2}
\end{figure}

Given these system and demand parameters, we first obtained $\overline{\omega}_d$ and $\underline{\omega}_d$ using 
$\left(\overline{\omega}^\prime_d, \underline{\omega}^\prime_d, \overline{r}_d, \underline{r}_d\right)$ 
for each $d$ according to the procedure described in Section \ref{sec:howtoomega}, 
which led to $M=64$. For instance, Fig. \ref{fig:2} depicts $\left(\overline{\omega}^\prime_1, \underline{\omega}^\prime_1\right)$ and $\left(\overline{\omega}_1, \underline{\omega}_1\right)$.

Subsequently, we solved (\ref{eq:prob_our_LP}) to determine $\left(\alpha,\beta_{1:M}\right)$ using the cutting-plane method, which required eight iterations and yielded the optimal value, i.e., the worst-case total generation 
cost, of $8.63\times10^4$. The lowest possible total generation cost equals $7.58\times10^4$.

\begin{figure}[t!]
    \centering
       \includegraphics[width=8.5cm]
       {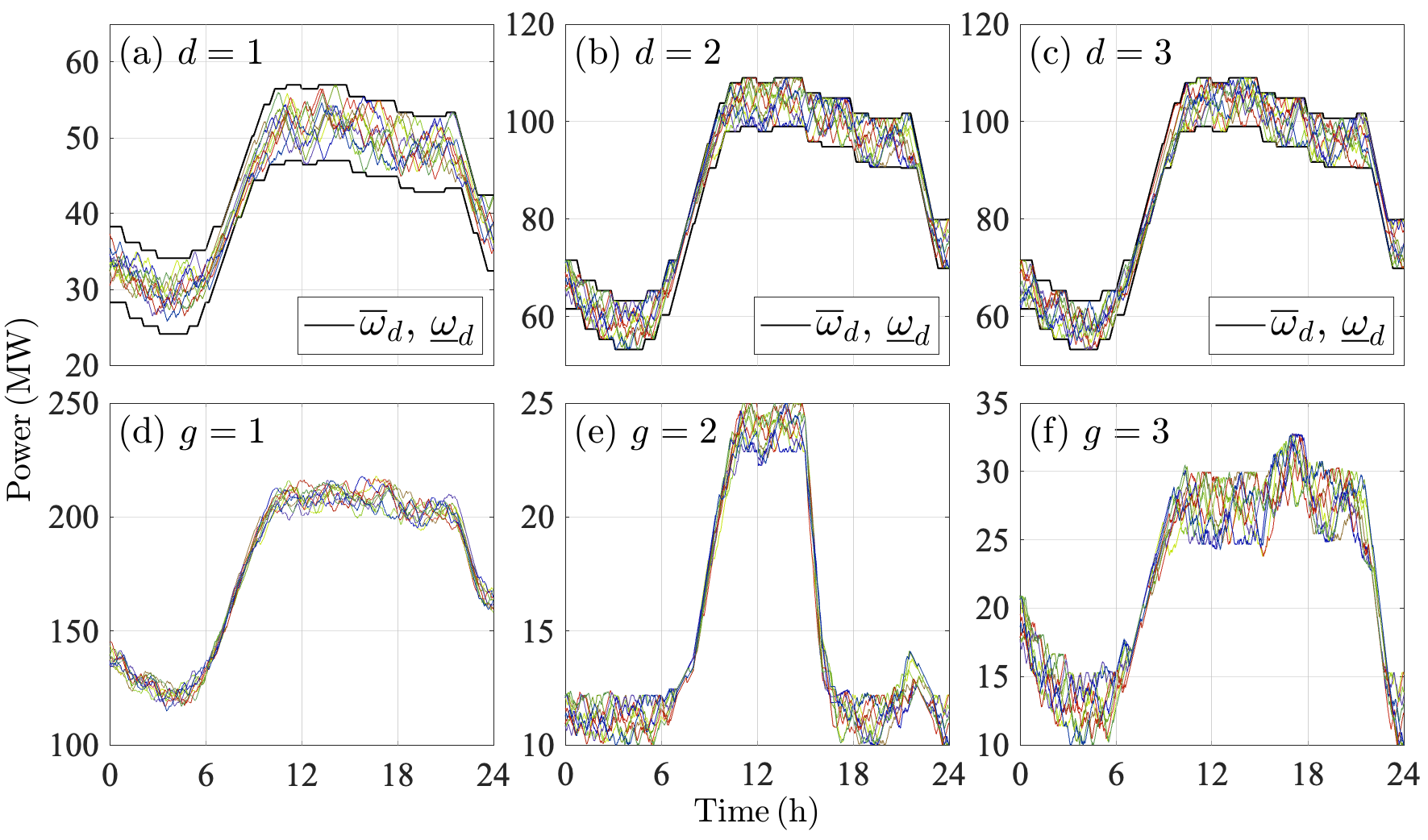}
    \caption{(a)--(c) The test demand trajectories. 
    (d)--(f) The power output trajectories obtained using our decision rule for the test demand trajectories.}
\label{fig:3}
\end{figure}

Our decision rule can accommodate any demand trajectory in $\Omega$. For example, we randomly generated ten test trajectories from $\Omega$ in the form of an affine function with $2401$ equally spaced breakpoints over $\mathcal T$. These demand trajectories and the corresponding power output trajectories produced by our decision rule are shown in Fig. \ref{fig:3}, the average total cost of which is equal to $8.15\times10^4$. Notably, each point of the power output trajectories were obtained non-anticipatively, using only the current demand values.

\begin{figure}[t!]
    \centering
       \includegraphics[width=8.5cm]
       {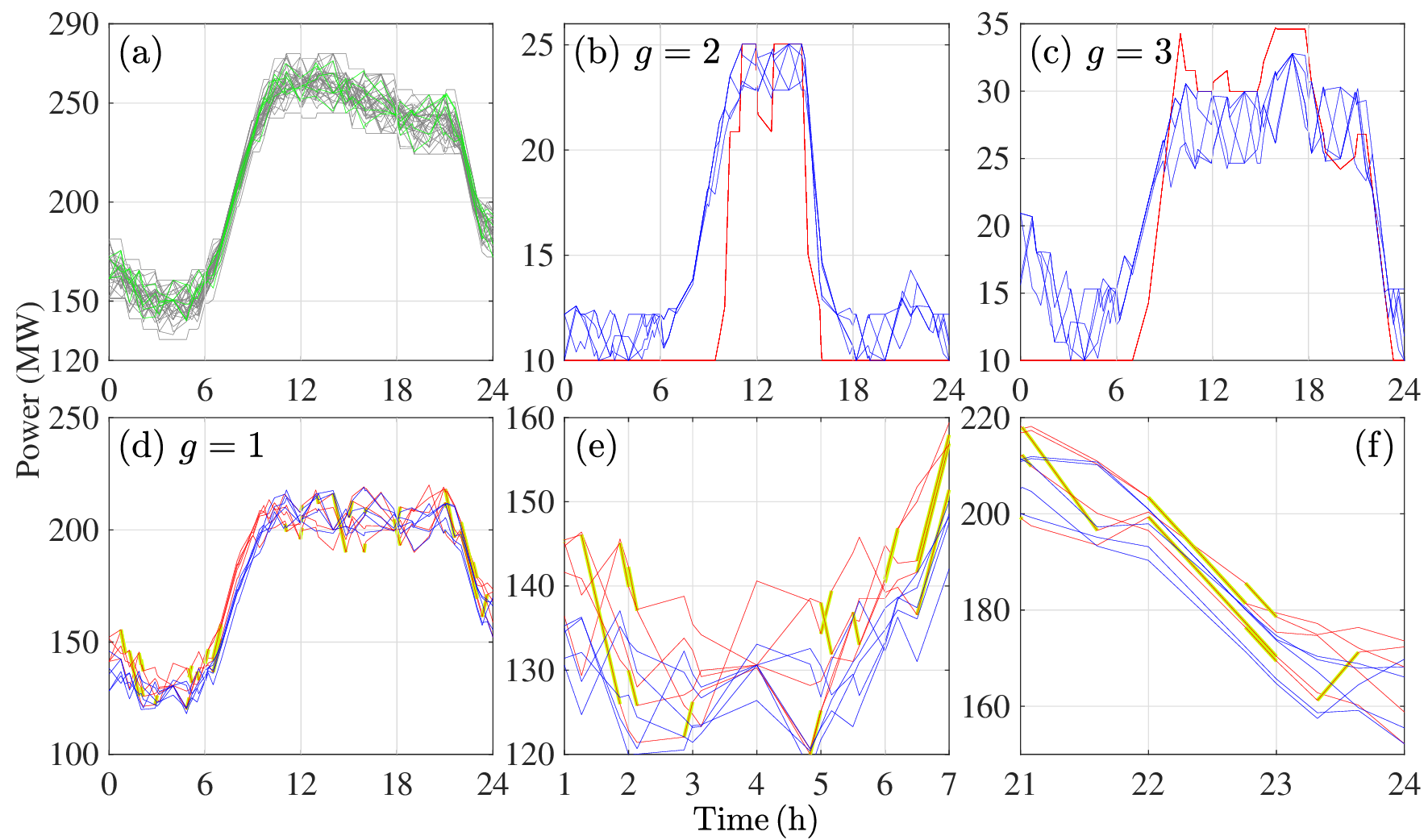}
    \caption{(a) The total demand computed with the trajectories in $\hat\Omega^{\rm sc}$ (grey) and the test trajectories from $\Omega$ (green). 
    (b)--(d) The power output trajectories obtained using $\hat{x}^{\rm sc}$ (red) and our decision rule (blue) for the test demand trajectories. (e), (f) Enlarged views of (d). 
    The yellow segments in (d)--(f) indicate the ramp-rate limit violations of generator 1, caused by $\hat{x}^{\rm sc}$.}
\label{fig:4}
\end{figure}

To emphasize the robustness of our model by way of contrast, we also solved a scenario-based counterpart. This model was obtained by replacing $\Omega$ in (\ref{eq:prob_our}) with $\hat{\Omega}^{\rm sc}:=
\{
\overline{\omega},\underline{\omega},\hat{\xi}_1,\ldots,\hat{\xi}_{30}
\}$. 
Here, $\hat{\xi}_1,\ldots,\hat{\xi}_{30}$ denote 30 piecewise affine functions randomly selected from $\Omega$, the breakpoint sets of which are ${\mathcal T}^{\rm b}$. Let $\hat{x}^{\rm sc}$ denote the resulting decision rule. The time-invariant coefficients of $\hat{x}^{\rm sc}$ for generators 2 and 3 are zero, implying that only generator 1, with the lowest marginal generation cost, adapts to demand uncertainty. 
Using $\hat{x}^{\rm sc}$, we obtained the power output trajectories of the generators non-anticipatively for another five test trajectories from $\Omega$. However, all the five trajectories of generator 1 violate the ramp rate constraints. 
Fig. \ref{fig:4} shows the trajectories of total demand for $\hat{\Omega}^{\rm sc}$ and the five test trajectories, in addition to the power output trajectories obtained using $\hat{x}^{\rm sc}$ for the test trajectories. The figure illustrates the power output trajectories obtained using our decision rule as well, which satisfy all the operational constraints of (\ref{eq:prob_our}) by design.

\section{Conclusions}\label{sec:conclusions}
We proposed a robust continuous-time generation scheduling model under demand uncertainty using a decision rule approach. The significance lies in its non-anticipative nature, distinguishing it from most existing continuous-time generation scheduling models. Future research directions include developing a systematic method for constructing the set of demand trajectories, which is presumed to be given in this study, and extending the proposed model to incorporate the operating status of generators as decision variables.


\appendix 
\subsection{Proof of Proposition \ref{prop:cost}}
\label{app:prop:cost}
Since $C^\top \alpha_{d}$ is constant over $\mathcal T$ for any $d\in{\mathbb N}\left[D\right]$ and 
\[
\left\{\xi\left(t\right):\xi\in\Omega\right\}=
\left[ \underline\omega_1\left(t\right), \overline\omega_1\left(t\right) \right] \times \cdots \times \left[ \underline\omega_D\left(t\right),\overline\omega_D\left(t\right) \right]
\]
for any $t\in{\mathcal T}$, 
$f_\Omega\left(\hat{x}_{\alpha,\beta_{1:M}}\right)$ 
is equal to 
\[
\begin{aligned}
&\sup_{\xi\in\Omega} \int_{{\mathcal T}} C^\top\alpha\xi\left(t\right) + C^\top \beta\left(t\right)dt= \int_{{\mathcal T}} C^\top \beta\left(t\right)dt\\ 
&+\sum_{d=1}^D\sup\bigg\{ 
\int_{\mathcal T}
C^\top\alpha_d \overline{\omega}_d\left(t\right)dt, \int_{\mathcal T}
C^\top\alpha_d\underline{\omega}_d\left(t\right)dt 
\bigg\}.
\end{aligned}
\]
As $\overline\omega$, $\underline{\omega}$, and $\beta$ are piecewise affine, we have 
\[
\begin{aligned}
&\int_{{\mathcal T}} C^\top \beta\left(t\right)dt=
\frac{1}{2}
\sum_{j=1}^{M-1}\Delta^{\rm b}_j
C^\top(\beta_{j+1}+\beta_j),\\
&\int_{\mathcal T}
C^\top\alpha_d \overline{\omega}_d\left(t\right)dt=
\frac{1}{2}\sum_{j=1}^{M-1}\Delta^{\rm b}_j
\overline{\omega}^{\rm s}_{d,j}
C^\top\alpha_d\quad \forall d\in{\mathbb N}\left[D\right],\\
&\int_{\mathcal T}
C^\top\alpha_d \underline{\omega}_d\left(t\right)dt=
\frac{1}{2}\sum_{j=1}^{M-1}\Delta^{\rm b}_j \underline{\omega}^{\rm s}_{d,j}
C^\top \alpha_{d}\quad \forall d\in{\mathbb N}\left[D\right].
\end{aligned}
\]
Hence, the statement holds.

\subsection{Proof of Proposition \ref{prop:ineq}}\label{app:prop:ineq}
We first prove the sufficiency. Condition (\ref{eq:prob_our_con_ineq}) is equivalently rewritten as 
\[
\left(A\alpha+B\right) \xi\left(t\right) + A\beta\left(t\right) \le a\quad 
\forall \left(\xi,t\right)\in\Omega\times{\mathcal T},
\] 
implying that 
\[ 
\left(A\alpha+B\right) \zeta + A\beta_j \le a
\quad 
\forall \zeta \in 
\left\{\xi\left(T^{\rm b}_j\right):\xi\in\Omega\right\}
\]
for any $j\in{\mathbb N}\left[M\right]$. Meanwhile, $(\underline{\omega}\left(T^{\rm b}_j\right) + \sigma_l \circ H_{j})\in{\mathbb R}^D$ for any $\left(j,l\right)\in{\mathbb N}\left[M\right]\times{\mathbb N}\left[2^D\right]$ corresponds to a vertex of 
\[ 
{\mathcal B}_j:=\left[\underline{\omega}\left(T^{\rm b}_j\right),\overline{\omega}\left(T^{\rm b}_j\right)\right]=\left\{\xi\left(T^{\rm b}_j\right):\xi\in\Omega\right\}.
\]
More precisely, we have 
\begin{equation}\label{eq:pf2_int}
\left\{\underline\omega\left(T^{\rm b}_j\right) + \sigma_l\circ
H_{j}: l\in{\mathbb N}\left[2^D\right]\right\}={\mathcal V}\left({\mathcal B}_j\right)\subset{\mathcal B}_j
\end{equation}
for any $j\in{\mathbb N}\left[M\right]$, indicating that the sufficiency holds. 

Subsequently, we prove the necessity. Due to the equality in (\ref{eq:pf2_int}),  
(\ref{eq:alphabeta_both_cap_lines}) implies that
\begin{subequations}
\begin{align}
&\left(A\alpha + B \right) \zeta + A\beta_j \le a\quad\forall \zeta\in{\mathcal B}_{j}, \label{eq:pf2_nec1}\\
&\left(A\alpha + B \right) \zeta^\prime + A\beta_{j+1} \le a\quad\forall \zeta^\prime\in{\mathcal B}_{j+1}\label{eq:pf2_nec2}
\end{align} 
\end{subequations}
for any $j\in{\mathbb N}\left[M-1\right]$. By taking a weighted sum of (\ref{eq:pf2_nec1}) and (\ref{eq:pf2_nec2}), we obtain 
\[
\begin{aligned}
&\left(A\alpha + B \right)\left(
\left(1-\delta\right) \zeta + \delta\zeta^\prime
\right)+A
\left[
\left(1-\delta\right)\beta_j + 
\delta\beta_{j+1} 
\right]
\le a\\
&\quad\forall \left(\zeta,\zeta^\prime,\delta\right)\in{\mathcal B}_j\times{\mathcal B}_{j+1}\times\left[0,1\right].
\end{aligned}
\]
This suggests that 
\[
\begin{aligned}
&\left(A\alpha
+ B \right)
\zeta + A
\left[
\left(1-\gamma\left(t\right)\right)\beta_{k\left(t\right)} + 
\gamma\left(t\right)\beta_{{k\left(t\right)}+1} 
\right]
\le a\\
&\quad\forall 
\zeta \in \left(1-\gamma\left(t\right)\right){\mathcal B}_{k\left(t\right)} + \gamma\left(t\right){\mathcal B}_{{k\left(t\right)}+1}
\end{aligned}
\]
for any $t\in{\mathcal T}$. As we have
\[
\begin{aligned}
&\left(1-\gamma\left(t\right)\right){\mathcal B}_{k\left(t\right)} + \gamma\left(t\right){\mathcal B}_{k\left(t\right)+1}\\
&=
\left\{
\left(1-\gamma\left(t\right)\right)\zeta + \gamma\left(t\right)\zeta^\prime:
\zeta\in{\mathcal B}_{k\left(t\right)}, 
\zeta^\prime\in{\mathcal B}_{k\left(t\right)+1}
\right\}
\\
&=
\left[\underline{\omega}\left(t\right),\overline{\omega}\left(t\right)\right]=\left\{
\xi\left(t\right): \xi\in\Omega
\right\}
\end{aligned}
\]
for any $t\in{\mathcal T}$, (\ref{eq:prob_our_con_ineq}) follows. This completes the proof.

\subsection{Proof of Proposition \ref{prop:ramp}}
\label{app:prop:ramp}
We first prove the sufficiency. For any $j\in{\mathbb N}\left[M-1\right]$, we have
\begin{equation}\label{eq:g_j}
g_j = \frac{\beta\left(t_2\right) - \beta\left(t_1\right)}{t_2-t_1}
\quad \forall\left(t_1, t_2\right)\in{\mathcal T}^{\rm o}_j 
\end{equation}
where
\[
{\mathcal T}^{\rm o}_j:=\left\{\left(t_1,t_2\right)\in{\mathcal T}^2: 
T^{\rm b}_j\le t_1< t_2\le T^{\rm b}_{j+1}\right\}
\]
denotes the set of all strictly ordered pairs of time points in the time interval $[T^{\rm b}_j, T^{\rm b}_{j+1}]$. Hence, (\ref{eq:prob_our_con_ramp}) implies that for any $j\in{\mathbb N}\left[M-1\right]$, 
\begin{equation}\label{eq:foaifjeo}
\underline{R}\le \alpha \frac{\xi\left(t_2\right)-\xi\left(t_1\right)}{t_2 - t_1} + g_j
\le \overline{R} \quad \forall \left(\xi,t_1,t_2\right)\in\Omega\times{\mathcal T}^{\rm o}_j. 
\end{equation}
For any $\left(d,j,\sigma\right)\in{\mathbb N}\left[D\right]\times{\mathbb N}\left[M-1\right]\times\left\{0,1\right\}$, there exists some $\left(\xi_d,t_1,t_2\right)\in\Omega\times{\mathcal T}^{\rm o}_j$ such that 
\[
R\left(\xi_d,t_1,t_2\right) =
\underline{r}^\prime_{j,d} + \sigma \left(\overline{r}^\prime_{j,d}-\underline{r}^\prime_{j,d}\right).
\] 
This equation trivially holds if $\overline{\omega}_d (T^{\rm b}_j)=\underline{\omega}_{d} (T^{\rm b}_j)$ 
and 
$\overline{\omega}_d (T^{\rm b}_{j+1})=\underline{\omega}_{d} (T^{\rm b}_{j+1})$, implying $\overline{r}^\prime_{j,d}=\underline{r}^\prime_{j,d}$. Otherwise, we can find such $\left(\xi_d, t_1, t_2\right)$ by considering 
any $t_1>T^{\rm b}_j$ and $t_2<T^{\rm b}_{j+1}$ that are arbitrarily close.  
Thus, (\ref{eq:alphabeta_ramp}) follows.

Subsequently, we prove the necessity. If (\ref{eq:alphabeta_ramp}) holds, then for any $j\in{\mathbb N}\left[M-1\right]$, we have  
\[
\underline{R} \le \alpha r + g_j \le \overline{R}
\quad
\forall r\in\left[\underline{r}^\prime_j,\overline{r}^\prime_j\right]
\]
as $\left(
\underline{r}^\prime_j + \sigma_l\circ (\overline{r}^\prime_j-\underline{r}^\prime_j)\right)$ for any $l\in{\mathbb N}\left[2^D\right]$ corresponds to a vertex of $\left[\underline{r}^\prime_j,\overline{r}^\prime_j\right]$, or more precisely, 
\[
\left\{\underline{r}^\prime_j
+\sigma_l\circ\left(\overline{r}^\prime_j-\underline{r}^\prime_j\right): l\in{\mathbb N}\left[2^D\right]\right\}
= {\mathcal V}\left(\left[\underline{r}^\prime_j,\overline{r}^\prime_j\right]\right).
\]
For any $j\in{\mathbb N}\left[M-1\right]$, since we have (\ref{eq:g_j}) and
\[
\frac{\xi\left(t_2\right) - \xi\left(t_1\right)}{t_2 - t_1}
\in\left[\underline{r}^\prime_j,\overline{r}^\prime_j\right]
\quad
\forall \left(\xi,t_1,t_2\right)\in\Omega\times{\mathcal T}^{\rm o}_j
\]
by definition, (\ref{eq:foaifjeo}) follows. This implies that 
\begin{equation}\label{eq:faefwefo}
\underline{R}\left(t_2-t_1\right)
\le
\hat{x}\left(\xi,t_2\right)-\hat{x}\left(\xi,t_1\right)
\le
\overline{R}\left(t_2-t_1\right)
\end{equation}
for any $\left(\xi,t_1,t_2\right)\in\Omega\times{\mathcal T}^{\rm o}$ such that $k\left(t_1\right)= k\left(t_2\right)$. 
Meanwhile, for any $\left(\xi,t_1,t_2\right)\in\Omega\times{\mathcal T}^{\rm o}$ such that $k\left(t_1\right)< k\left(t_2\right)$, we obtain from (\ref{eq:foaifjeo}) that  
\[
\begin{aligned}
\hat{x}\left(\xi,{t_2}\right) -\hat{x}\left(\xi, T_{k\left(t_2\right)}\right)
&\le \overline{R}\left(t_2 - T_{k\left(t_2\right)}\right),\\
\hat{x}\left(\xi,T_{k\left(t_2\right)}\right) -\hat{x}\left(\xi, T_{k\left(t_2\right)-1}\right)
&\le \overline{R}\left(T_{k\left(t_2\right)} - T_{k\left(t_2\right)-1}\right),\\ 
&\cdots\\
\hat{x}\left(\xi,T_{k\left(t_1\right)+2}\right) -\hat{x}\left(\xi, T_{k\left(t_1\right)+1}\right)&\le\overline{R}\left(T_{k\left(t_1\right)+2} - T_{k\left(t_1\right)+1}\right), \\
\hat{x}\left(\xi,T_{k\left(t_1\right)+1}\right) -\hat{x}\left(\xi, t_1\right)
&\le \overline{R}\left(T_{k\left(t_1\right)+1} - t_1\right)
\end{aligned}
\]
where the first inequality trivially holds when $t_2=T_{k\left(t_2\right)}$. Adding these $k\left(t_2\right)-k\left(t_1\right)+1$ inequalities, we have
\[
\hat{x}\left(\xi,t_2\right) - \hat{x}\left(\xi,t_1\right)\le
\overline{R}\left(t_2-t_1\right).
\] 
Similarly, we can derive the inequality
\[
\underline{R}\left(t_2-t_1\right)\le
\hat{x}\left(\xi,t_2\right) - \hat{x}\left(\xi,t_1\right). 
\]
Thus, (\ref{eq:faefwefo}) holds for any $\left(\xi,t_1,t_2\right)\in\Omega\times{\mathcal T}^{\rm o}$, which is equivalently expressed by (\ref{eq:prob_our_con_ramp}). The proof is complete.

\subsection{Proof of Proposition \ref{prop:bal}}\label{app:prop:bal}
As the necessity is obvious, we prove only the sufficiency. Constraint (\ref{eq:prob_our_con_bal}) is rewritten as 
\begin{equation}\label{eq:lffoissi}
\begin{aligned}
&\left({1}^\top_G\alpha-{1}^\top_D\right)\xi\left(t\right) 
+ \left(1-\gamma\left(t\right)\right){1}^\top_G\beta_{k\left(t\right)} \\
&\quad + \gamma\left(t\right){1}^\top_G\beta_{k\left(t\right)+1}= 0 \quad
\forall \left(\xi,t\right)\in\Omega\times{\mathcal T}. 
\end{aligned}
\end{equation}
For any $d\in{\mathbb N}\left[D\right]$, there exists some $t\in{\mathcal T}$ such that $\left\{\xi_d\left(t\right):\xi_d\in\Omega_d\right\}=\left[\underline{\omega}_d\left(t\right),\overline{\omega}_d\left(t\right)\right]$ is non-degenerate as $\Xi_d$ is non-singleton. Thus, there exists some non-zero $b\in{\mathbb R}$ such that $b\left({1}^\top_G\alpha - {1}^\top_D\right){e}_d=0$, where ${e}_d$ denotes the $d$th column of the $D \times D$ identity matrix, implying ${1}^\top_G\alpha - {1}^\top_D=0$. 
Using this equation, (\ref{eq:lffoissi}) is rewritten as 
\begin{equation}\label{eq:flafif}
\left(1-\gamma\left(t\right)\right){1}^\top_G\beta_{k\left(t\right)} + \gamma\left(t\right){1}^\top_G\beta_{k\left(t\right)+1}= 0 \quad\forall t\in{\mathcal T}. 
\end{equation} 
Thus, we have 
\[
\left(\gamma\left(t_2\right)-\gamma\left(t_1\right)
\right)
\left( {1}^\top_G\beta_{k\left(t_1\right)} - 
{1}^\top_G\beta_{k\left(t_1\right)+1}
\right) = 0
\] 
for any $\left(t_1,t_2\right)\in{\mathcal T}^{\rm o}$ such that $k\left(t_1\right)=k\left(t_2\right)$, which implies 
\[
{1}^\top_G \beta_{k\left(t_1\right)} = {1}^\top_G\beta_{k\left(t_1\right)+1}\]
as $\gamma\left(t_2\right)\neq\gamma\left(t_1\right)$. 
Due to (\ref{eq:flafif}), we observe 
\[
{1}^\top_G \beta_{k\left(t_1\right)} = {1}^\top_G \beta_{k\left(t_2\right)} = 0.
\]
This completes the proof.

\bibliographystyle{IEEEtran}
\bibliography{reference}

\end{document}